\documentclass[12pt,a4paper]{article}
\usepackage{color}
\usepackage{t1enc}
\usepackage[dvips,final]{graphicx} 
\usepackage{tabularx} 
\usepackage{amssymb,amsfonts,amsmath,amscd,euscript}
\usepackage{graphics}
\usepackage{psfrag}
\usepackage{caption}
\usepackage{amsthm}

\def\softd{{\leavevmode\setbox1=\hbox{d}%
\hbox to 1.05\wd1{d\kern-0.4ex{\char039}\hss}}}
\def\softt{{\leavevmode\setbox1=\hbox{t}%
\hbox to \wd1{t\kern-0.6ex{\char039}\hss}}}
\def\softl{l\kern-0.45ex\raise0.1ex\hbox{'}\kern-0.10ex}
\def\softL{L\kern-0.8ex\raise0.1ex\hbox{'}\kern0.1ex}

\newcommand\curl{\text{\rm curl\,}}
\newcommand\Div{{\text{\rm div\,}}}
\newcommand\grad{\text{\rm grad\,}}
\newcommand\Ham{\text{\rm Ham\,}}

\begin{document}
\title{Vector calculus in two-dimensional space}
\author{Mari\'an Fecko\thanks{e-mail: {\texttt{Marian.Fecko@fmph.uniba.sk}}} \\
\small{{Department of Theoretical Physics}}                                \\
\small{{Faculty of Mathematics, Physics and Informatics}}                  \\
\small{{Comenius University in Bratislava}}}
\date{}
\maketitle
\small{\emph{Vector calculus in three-dimensional space is ubiquitous in
             applications of mathe\-matics in physics and engineering.
Its two-dimensional version is, however, quite rare.
Here we try to provide a pedagogical account of the subject.
It is based on the logic of theory of differential forms.
For readers not familiar with the latter, the results are presented in detail
in the standard language and notations.
             }
\vskip .5 cm

\vglue -1 cm
\tableofcontents

\newpage

\section{Introduction}
\label{intro}

We live in three-dimensional space.
Therefore also mathematics connected with the three-dimensional space turns out to be especially important
as well as developed.

In particular, it turns out that \emph{vector calculus} ideally suits for neat and clear mathematical description
of a number of important phenomena around us
(diffusion, heat flow, electricity, magnetism, gravitational field, ...; see e.g.
\cite{Fleisch}
\cite{Marsden},
\cite{Schey},
\cite{Spain},
 etc.).
Useful concepts are introduced in it (like flux of a vector field for a surface)
and useful and non-trivial results (as, for example, Gauss' theorem) are derived for them.

Two dimensional version of the vector calculus is, however, considerably less known.
It should play the same role for 2D phenomena as the standard vector calculus does in 3D.
Nevertheless, what and how exactly is to be changed in mathe\-ma\-tical relations
while switching from 3D to 2D may not always be quite evident.

In this text we therefore concentrate on this very question:
How actually \emph{the} analogy of the three dimensional vector calculus,
applicable in \emph{two} dimensions, looks like?

Several approaches might exist, see e.g. \cite{Khan} or \cite{Olver}.
Here wee try to be led by logic of the theory of \emph{differential forms}.
For reader unfamiliar with the theory it might look as an overkill.
It is, however, well-known (see e.g.
\cite{BambergSternberg},
\cite{CrampinPirani},
\cite{Schutz} or \$8.5 in
\cite{Fecko}),
that this very point of view turns out to be extremely instructive and efficient in~\emph{three}-dimensional case:
For example, all relevant \emph{differential operators} do appear automatically
as well as relations between them (differently looking curl, grad and div are just masked versions
of universal \emph{exterior} derivative of forms).
This is also the case for all \emph{integral theorems} in which the operators occur
(as particular cases of each time the same \emph{Stokes} theorem for forms).

In this text we can see that differential forms approach is definitely not a bad choice also here.

If forms in general (and so also vector calculus in the language of forms in particular)
is virtually unknown to the reader, derivations might be technically unclear, but:

1. \emph{Results} will be clear, since they will be also presented in \emph{elementary} language,
   so that just to understand \emph{how it all falls out}, forms are not needed.

2. In addition, rough idea perhaps still might be caught,
   and consequently the forms could captivate the reader, resulting in decision
   to learn more about them in the future.

\section{How it works in three dimensions}
\label{akov3d}

In vector calculus one often hears phrases like ``integrate vector field over a surface''
(resulting in its \emph{flux} for the surface).
In integration theory, however, we learn (see e.g.
\cite{BambergSternberg},
\cite{CrampinPirani},
\cite{Schutz} or
Chapter 7 in \cite{Fecko}), that

- differential \emph{forms} alone are integrated

- $p$-form is integrated over a $p$-dimensional domain.

\noindent
So if something is integrated

- along a curve, it is necessarily a 1-form,

- over a surface, it is necessarily a 2-form,

- over a volume, it is necessarily a 3-form, respectively.

\noindent
How should one understand the fact, then, that under surface integrals, where 2-form is to be standing,
one sees expression of the structure $\mathbf A \cdot d\mathbf S$,
i.e. not to see \emph{vector field} $\mathbf A$ there is like even deny what everything someone does?
Well, the answer is that what we \emph{actually} see there is \emph{the whole} $\mathbf A \cdot d\mathbf S$
and this in fact \emph{is} a 2-form.
It is, however, \emph{uniquely parametrized} in terms of the vector field $\mathbf A$.

In total, parametrizations of all types of forms in $E^3$ (three dimensional Euclidean space) look as follows
(see e.g. section 8.5 in \cite{Fecko} mentioned above):

There are 0-forms, 1-forms, 2-forms and 3-forms (this is just because of three dimensions of the space)
and we have, at each of these degrees, operation of \emph{exterior derivative} $d$,
raising the degree in one unit and such that its square vanishes (this is so in any dimension).
One can display the situation in terms of the diagram
\begin{equation} \label{DeRhamE3}
          \Omega^0  \overset d \to  \Omega^1    \overset d \to    \Omega^2  \overset d \to   \Omega^3
          \hskip 2cm
          dd=0
\end{equation}
(it is known as the \emph{de Rham complex}; the space of $p$-forms in $E^3$ is denoted as $\Omega^p$).
Since there is also (usual Euclidean) \emph{metric tensor} on $E^3$, we have in addition
the \emph{Hodge operator} on $p$-forms. It is, in general, a canonical linear isomorphism
(denoted by star) of the linear space of $p$-forms and $(n-p)$-forms,
which squares to plus or minus identity.
(It is just plus on each degree on $E^3$. Here $n$ is dimension of the space under consideration.)

So in \emph{three}-dimensional space it provides the following two canonical isomorphisms:
\begin{equation} \label{HodgeE3}
          \Omega^0  \overset * \leftrightarrow  \Omega^3
          \hskip 1.5cm
          \Omega^1  \overset * \leftrightarrow  \Omega^2
          \hskip 1.5cm
          *^{-1} = *
\end{equation}
They enable us to ``identify'' 2-forms with 1-forms as well as 3-forms with 0-forms.
(If we need somewhere a 2-form $\beta$, we can \emph{parametrize} it in terms of a 1-form $\alpha$,
 i.e. write it as $\beta = *\alpha$ where $\alpha$ is unique 1-form: $\alpha = *\beta$.)
So ``actually'' we need just 0-forms (i.e. \emph{functions}; we denote as $\mathcal F$ the linear space of functions)
and 1-forms (the rest may be parametrized in terms of the two).

The metric tensor provides another two important canonical isomorphisms, namely the (mutually inverse)
operations of the \emph{raising and lowering of indices}
(they are denoted as $\sharp$ and $\flat$, inspiration originating in musical notation).
They turn 1-forms into vector fields ($\sharp$) and vice versa ($\flat$).
If the linear space of vector fields on $E^3$ is denoted as $\mathfrak X$, we have
\begin{equation} \label{flatsharp}
          \Omega^1  \underset \flat {\overset \sharp \rightleftarrows}   \mathfrak X
          \hskip 2cm
          \sharp^{-1} = \flat
          \hskip 1cm
          \flat^{-1} = \sharp
\end{equation}
This is the reason we can also forget about 1-forms and manage with
\footnote{In this respect functions and vector fields resemble concept of  \emph{essential} amino-acids in nutrition science:
          Although there are roughly 20 amino-acids in total, it is actually enough to take care of the essential ones
          (roughly half of them), our body is then able to produce the remaining ones from them.
          In three dimensional Euclidean space our body is able to produce \emph{any} differential form
          from functions and vector fields.}
just \emph{functions and vector fields}! So differential forms in $E^3$ become \emph{as if} completely out of the game!
We can concisely express the situation in terms of the commutative diagram
\begin{equation}
\begin{CD} \label{diagram3D1}
          \Omega^0          @>d>>   \Omega^1       @>d>>    \Omega^2        @>d>>    \Omega^3    \\
           @V{\text id}VV            @VV{\sharp}V            @VV{\sharp *}V           @VV*V      \\
          \mathcal F    @>>a_0>     \mathfrak X     @>>a_1>     \mathfrak X    @>>a_2>    \mathcal F  \\
\end{CD}
\end{equation}
Vertical arrows denote the corresponding canonical isomorphisms of various degrees of forms
(the upper line) onto scalar fields $\mathcal F$ and vector fields $\mathfrak X$.
When the arrows are reversed, we get equivalent diagram
\begin{equation}
\begin{CD} \label{diagram3D2}
          \Omega^0          @>d>>   \Omega^1       @>d>>    \Omega^2        @>d>>    \Omega^3    \\
           @A{\text id}AA            @AA{\flat}A            @AA{*\flat}A          @AA*A      \\
          \mathcal F    @>>a_0>     \mathfrak X     @>>a_1>     \mathfrak X    @>>a_2>    \mathcal F  \\
\end{CD}
\end{equation}
When the (upward directed) arrows are applied onto scalar and vector fields (on some $f$ and $\mathbf A$), we get
exactly those standard expressions we used to see under integral sign (see section 8.5 in \cite{Fecko}):
\begin{equation} \label{formsstandardne}
          f \in \Omega^0
          \hskip 1cm
          \mathbf A \cdot d\mathbf r \in \Omega^1
          \hskip 1cm
          \mathbf A \cdot d\mathbf S \in \Omega^2
          \hskip 1cm
          fdV \in \Omega^3
\end{equation}
The arrows $a_0,a_1,a_2$ in the bottom line of both diagrams are \emph{effective} operations
(composition of corresponding three arrows - upward, to the right and downward)
on objects sitting in the bottom line, which in a sense ``substitute'' operation of the exterior derivative
(arrows $d$), which ``really'' acts on forms in the upper line.
From the diagrams we can read off that
\begin{equation} \label{a0a1a2}
          a_0 = \sharp d
          \hskip 1cm
          a_1 = \sharp * d \flat
          \hskip 1cm
          a_2 = *d *\flat
\end{equation}
All of them are clearly \emph{first} order \emph{differential} operators, since they contain one $d$.
Computation in Cartesian coordinates in $E^3$ reveals that they are exactly notorious operators
grad, curl and div so that the diagram (\ref{diagram3D1}) actually reads
\begin{equation}
\begin{CD} \label{diagram3D3}
          \Omega^0          @>d>>   \Omega^1       @>d>>    \Omega^2        @>d>>    \Omega^3    \\
           @V{\text id}VV            @VV{\sharp}V            @VV{\sharp *}V           @VV*V      \\
          \mathcal F    @>>\grad>     \mathfrak X     @>>\curl>     \mathfrak X    @>>\Div>    \mathcal F  \\
\end{CD}
\end{equation}
From the fact that $dd=0$ (see (\ref{DeRhamE3})) and commutativity of the diagram we immediately
get well-known identities
\begin{equation} \label{znameidentity}
          \curl \grad = 0
          \hskip .7cm
          \Div \curl = 0
\end{equation}
(upstairs, composition of neighboring arrows results in zero operator,
 so the same must be true downstairs).

From the three neighboring squares of the diagram (\ref{diagram3D3}) and the parametrization (\ref{formsstandardne})
one can assemble the following three useful differential relations:
\begin{equation} \label{dnaformsstandardne}
          df = \grad f \cdot d\mathbf r
          \hskip .7cm
          d(\mathbf A \cdot d\mathbf r)= (\curl \mathbf A ) \cdot d\mathbf S
          \hskip .7cm
          d(\mathbf A \cdot d\mathbf S)= (\Div \mathbf A ) dV
\end{equation}
In order to obtain corresponding integral statements, we can make reference
to general \emph{Stokes theorem} (see
\cite{BambergSternberg},
\cite{CrampinPirani},
\cite{Schutz} or Section 7.5 in
\cite{Fecko})
from integration theory of differential forms:
If $\alpha$ is a $p$-form, $D$ is a $(p+1)$-dimensional domain
and $\partial D$ its $p$-dimensional boundary, then
\begin{equation} \label{stokespreforms}
          \int_Dd\alpha = \int_{\partial D}\alpha
\end{equation}
For the three particular expressions $d\alpha$ from (\ref{dnaformsstandardne}) we get
the three essential integral theorems of the vector calculus:
\begin{align} 
              \label{gradient}
              \text{\emph{gradient theorem}}&  &  \int_c \grad f \cdot d\mathbf r            &:= f(B) - f(A)                                    \\
              \label{stokes}
              \text{\emph{Stokes   theorem}}&  &  \int_S (\curl \mathbf A ) \cdot d\mathbf S &:= \oint_{\partial S} \mathbf A \cdot d\mathbf r  \\
              \label{gauss}
              \text{\emph{Gauss's  theorem}}&  &  \int_V (\Div \mathbf A ) dV                &:= \oint_{\partial V} \mathbf A \cdot d\mathbf S
\end{align}

\section{How it works in two dimensions}
\label{akov2d}

In two-dimensional case, the de Rham complex, i.e. the analogue of (\ref{DeRhamE3}), just simplifies to
\begin{equation} \label{DeRhamE2}
          \Omega^0  \overset d \to  \Omega^1    \overset d \to    \Omega^2
          \hskip 2cm
          dd=0
\end{equation}
Key difference is in action of the Hodge star:
\begin{equation} \label{HodgeE2}
          \Omega^0  \overset * \leftrightarrow  \Omega^2
          \hskip 1.5cm
          \Omega^1  \overset * \leftrightarrow  \Omega^1
\end{equation}
The left expression is a natural analogue of the left expression in (\ref{HodgeE3}),
it identifies the edge degrees of forms; explicitly
\begin{equation} \label{HodgeE2okraje}
          * f = fdS
          \hskip 1cm
          * (fdS) = f
          \hskip 1.5cm
          dS \equiv dx \wedge dy = \ \text{area 2-form}
\end{equation}
This enables us to forget about 2-forms and express them in terms of 0-forms,
i.e. in terms of functions (scalar fields).

The right expression is, however, completely different:
While in (\ref{HodgeE3}) it is an isomorphism of \emph{two different} spaces,
in (\ref{HodgeE2}) it is an isomorphism of a \emph{single} space (the space of 1-forms) \emph{on itself}
(differing from the identity).
Explicitly, in usual Cartesian coordinates $(x,y)$ in the plane,
\begin{equation} \label{HodgeE2nabaze}
          * dx = dy
          \hskip 1.5cm
          * dy = -dx
\end{equation}
This does not mean that we cannot forget about 1-forms. Rather, it means that in 2D
there are as many as \emph{two canonical} ways how 1-forms may be replaced by vector fields!
One of them is the same like it was in 3D, i.e. via the operation $\sharp$
(raising of index). It reads
\begin{equation}
\begin{CD} \label{prvyanalog}
          \Omega^0      @>d>>   \Omega^1    @>d>>    \Omega^2       \\
           @V{\text id}VV                @VV{\sharp}V               @VV*V          \\
          \mathcal F    @>>a_0>   \mathfrak X     @>>a_1>    \mathcal F     \\
\end{CD}
\end{equation}
When the arrows (directed upwards) are applied on scalar and vector fields (on $f$ and $\mathbf A$)
we again get standard expressions we used to see under integral sign:
\begin{equation} \label{formsstandardne2D}
          f \in \Omega^0
          \hskip 1cm
          \mathbf A \cdot d\mathbf r \equiv A_xdx + A_ydy \in \Omega^1
          \hskip 1cm
          fdS \in \Omega^2
\end{equation}
Another one is, however, available, now - we first apply the star operator
(after which we still remain within 1-forms)
and \emph{only then} we raise the index.
So we have \emph{two different} analogues of the diagram (\ref{diagram3D1}):
\begin{equation}
\begin{CD} \label{twoanalogy}
          \Omega^0      @>d>>   \Omega^1    @>d>>    \Omega^2       \\
           @V{\text id}VV                @VV{\sharp}V               @VV*V          \\
          \mathcal F    @>>a_0>   \mathfrak X     @>>a_1>    \mathcal F     \\
\end{CD}
\hskip 1.5cm
\begin{CD}
          \Omega^0      @>d>>   \Omega^1    @>d>>    \Omega^2       \\
           @V{\text id}VV                @VV{\sharp *}V               @VV*V          \\
          \mathcal F    @>>b_0>   \mathfrak X     @>>b_1>    \mathcal F     \\
\end{CD}
\end{equation}
We see that also in 2D-vector calculus we make do with \emph{vector and scalar} fields (bottom lines).
From the diagrams one can express explicitly the bottom (effective) arrows and get
\begin{equation} \label{a0a1b0b1}
          a_0 = \sharp d
          \hskip 1cm
          a_1 =  * d \flat
          \hskip 1cm
          b_0 = \sharp * d
          \hskip 1cm
          b_1 = *d *^{-1}\flat
\end{equation}
When the (differential) operators are computed in (Cartesian) coordinates $(x,y)$,
we learn that two of them are natural analogues of the situation in 3D (gradient and divergence),
the remaining two are ``new'', specific for 2D:
\begin{equation}
\begin{CD} \label{twoanalogypodrobne}
          \Omega^0      @>d>>   \Omega^1    @>d>>    \Omega^2       \\
           @V{\text id}VV                @VV{\sharp}V               @VV*V          \\
          \mathcal F    @>>{\grad}>   \mathfrak X     @>>\curl_3>    \mathcal F     \\
\end{CD}
\hskip 1.5cm
\begin{CD}
          \Omega^0      @>d>>   \Omega^1    @>d>>    \Omega^2       \\
           @V{\text id}VV                @VV{\sharp *}V               @VV*V          \\
          \mathcal F    @>>{\Ham}>   \mathfrak X     @>>{- \Div}>    \mathcal F     \\
\end{CD}
\end{equation}
(Notation $b_0=\Ham$ is justified in section \ref{hamiltonfields} and for $a_1=\curl_3$ in section \ref{tretidimension}.)
So, abstractly we have
\begin{equation} \label{a0a1b0b1naozaj}
          \grad  = \sharp d
          \hskip 1cm
          \curl_3 =  * d \flat
          \hskip 1cm
          \Ham = \sharp * d
          \hskip 1cm
          \Div = - *d *^{-1}\flat
\end{equation}
and their actions, in Cartesian coordinates $(x,y)$, read:
\begin{align} 
              \label{grad}
              \grad :&  &  f         &\mapsto    (\partial_x f, \partial_y f)       \\
              \label{a1}
              \curl_3   :&  &  (A_x,A_y) &\mapsto    (\partial_x A_y -\partial_y A_x)   \\
              \label{Ham}
              \Ham   :&  &  f         &\mapsto    (-\partial_y f,\partial_x f)       \\
              \label{div}
              \Div  :&  &  (A_x,A_y) &\mapsto    (\partial_x A_x + \partial_y A_y)
\end{align}
One can easily check vanishing of composition of the two arrows in bottom lines (in both diagrams):
\begin{equation}
       \label{a1gradDivHam}
       \curl_3 \circ \grad = 0
       \hskip 1cm
       \Div \circ \Ham = 0
\end{equation}
Indeed,
\begin{align} 
              \label{a1gradvypocet}
              \curl_3 \circ \grad  :&  &  f  &\mapsto (\partial_x f, \partial_y f)    &  &\mapsto (\partial_y\partial_x f - \partial_y \partial_x f)   & &= 0      \\
              \label{divhamvypocet}
              \Div \circ \Ham   :&  &  f  &\mapsto ( - \partial_y f, \partial_x f) &  &\mapsto (- \partial_x \partial_y f, \partial_y \partial_x f) & &= 0
\end{align}
This is the counterpart of (\ref{znameidentity}). Notice that there are also two identities, here.
However, not because we put to use $dd=0$ on different degrees in a single diagram,
but rather on a single degree in two different diagrams.

It is also worth noticing that when two arrows from different diagrams are composed,
we get (for both possible cases) another well-known operator, the \emph{Laplace} operator $\triangle$:
\begin{equation}
       \label{gradDivHama1}
        \Div \circ \grad = \triangle
       \hskip 1cm
        \curl_3 \circ \Ham = \triangle
\end{equation}
Indeed,
\begin{align} 
              \label{Divgradvypocet}
              \Div \circ \grad   :&  &  f  &\mapsto (\partial_x f, \partial_y f)    &  &\mapsto (\partial^2_x + \partial^2_y)f & &\equiv \triangle f  \\
              \label{a1hamvypocet}
              \curl_3 \circ \Ham :&  &  f  &\mapsto ( - \partial_y f, \partial_x f) &  &\mapsto (\partial^2_x + \partial^2_y)f & &\equiv \triangle f
\end{align}

\subsection{Hamiltonian fields}
\label{hamiltonfields}

Notation $\Ham f$ in (\ref{twoanalogypodrobne}) and (\ref{Ham}) is abbreviation for \emph{Hamiltonian field} generated by function~$f$.
It is key concept in theory of Hamiltonian systems (for more details see, e.g.,
\cite{CrampinPirani},
\cite{Schutz} or Chapter 14 in \cite{Fecko}).
How could it occur here?

Well, it turns out that our \emph{surface} form $\omega \equiv dS$ mentioned in (\ref{HodgeE2okraje})
happens to be, at the same time, \emph{symplectic} form
(it satisfies corresponding definition; in general it is to be ``closed and non-degenerate 2-form'').
And whenever symplectic form is available,
a rule for construction of Hamiltonian (vector) fields is available as well:
\begin{equation}
       \label{hamfield}
       f\mapsto \Ham f
       \hskip 1cm
       i_{\Ham f}\omega := - df
\end{equation}
Now, when we compute how (\ref{hamfield}) works here, what we get is exactly (\ref{Ham}).

This also sheds a new light onto the second result in (\ref{a1gradDivHam}).
The fact, that divergence of Hamiltonian field (generated by any function) vanishes
is equivalent to well-known \emph{Liouville theorem} from classical mechanics (see 14.3.7 in \cite{Fecko}).
It says that, for Hamiltonian systems, \emph{volume} (here area)
in \emph{phase} space (here our 2D space) is conserved w.r.t. time evolution.

\subsection{Auxiliary third dimension}
\label{tretidimension}

One can easily check that:

\noindent
1. If, from a function $f(x,y)$ in 2D plane, one constructs the following auxiliary vector field in 3D space
\begin{equation}
       \label{fieldpreham}
             \mathbf u = (0,0,-f(x,y))
\end{equation}
(try to think how it looks like), then its curl reads
\begin{equation}
       \label{rotpolapreham}
             \curl \mathbf u =  (-\partial_y f,\partial_x f,0) \equiv (\Ham f,0)
\end{equation}
2. If, from a vector field $\mathbf A$ in 2D plane, one constructs the following auxiliary vector field in 3D space
\begin{equation}
       \label{fieldprerot3}
             \boldsymbol {\mathcal A} = (\mathbf A,0) \equiv (A_x(x,y),A_y(x,y),0)
\end{equation}
(try to think how it looks like), then its curl reads
\begin{equation}
       \label{rotpolaprerot3}
             \curl \boldsymbol  {\mathcal A} =  (0,0,\partial_x A_y - \partial_y A_x) \equiv (0,0,\curl_3 \mathbf A)
\end{equation}
What these two simple computations do reveal?

First, that our two ``new'' operations in 2D vector calculus,
$a_1 \equiv \curl_3$ and $b_0 \equiv \Ham$,
may also be discovered as hidden at special places of outputs of specially chosen inputs
of standard 3D vector calculus.
The approach treated in this paper finds their form and properties in ``intrinsic'' way,
with no reference to auxiliary additional dimensions.

And, second, we see motivation for notation of the operation $\curl_3$.
Well, mainly for the index 3
(in particular, that $\curl_3 \mathbf A = (\curl \boldsymbol {\mathcal A})_3$).
Motivation for the notation \emph{curl} itself (which should be related to kind of \emph{rotation} of something)
is discussed already in standard 3D vector calculus
(it stems for example from \emph{vortex}-like motion of fluid flow;
see a few more words on this in Section \ref{potencial}).

\subsection{Integral identities and Green's theorem}
\label{2dintegralnevety}

Integral theorems for 2D-vector calculus may be again derived from universal Stokes theorem
(\ref{stokespreforms}) for forms. Corresponding technical procedure just mimics the one used
in the 3D case at the end of Section \ref{akov3d}.

Since there are altogether as many as \emph{four} basic squares (with upper arrow $d$) in two diagrams
in (\ref{twoanalogypodrobne}), four integral theorems are expected.

First, we assemble the corresponding four differential identities
(analogues of expressions (\ref{dnaformsstandardne})).
From the individual squares we get
\begin{align} 
              \label{lavestvorce}
              df                              &=  \grad f  \cdot d\mathbf r  &    df                           &= - * (\Ham f  \cdot d\mathbf r ) \\
              \label{pravestvorce}
              d (\mathbf A \cdot d\mathbf r)  &=    (\curl_3 \mathbf A)dS        &  d*(\mathbf A \cdot d\mathbf r) &= (\Div \mathbf A ) dS
\end{align}
Their integration and application of Stokes theorem (\ref{stokespreforms}) leads to
\begin{align} 
\label{gradientova1}
\int_C\grad f  \cdot d\mathbf r &= f(B) - f(A)                                    &  \int_C * (\Ham f  \cdot d\mathbf r )   &= f(A) - f(B)  \\
\label{greenova1}
\int_S(\curl_3 \mathbf A)dS         &= \oint_{\partial S} \mathbf A \cdot d\mathbf r  &  \int_S(\Div \mathbf A ) dS             &= \oint *(\mathbf A \cdot d\mathbf r)
\end{align}
Now, when we explicitly write down expressions under integral signs (with the help of (\ref{HodgeE2nabaze}) and (\ref{grad}) - (\ref{div})),
we discover that both identities in (\ref{gradientova1}) say the same thing
and, similarly, the same is true for both identities in (\ref{greenova1}).
So we get \emph{just two} mutually different integral identities.
When appropriate notations are used, they read as follows:
\begin{align} 
              \label{gradient2d}
              \text{\emph{gradient theorem}}& & \int_c (\partial_x f)dx + (\partial_y f)dy  &:= f(B) - f(A)                  \\
              \label{green}
              \text{\emph{Green's theorem}}&   & \int_S (\partial_x g - \partial_y f)dxdy    &:= \oint_{\partial S} fdx +gdy
\end{align}
The first one is a 2D version of the gradient theorem. This is well known and holds in any dimension.
It relates an integral along a curve $c$ to an ``integral'' over its boundary $\partial c$
(which reduces, here, to just two points, $A$, its beginning and $B$, its end).

The second one is Green's theorem.
This is well-known, again, and is specific for two dimensional space.
It relates an integral over a 2D-domain $S$ to an integral over its boundary
(being a closed circuit $\partial S$).

\subsection{Perpendicular vector field}
\label{kolmefield}

As it was already mentioned in Section \ref{akov2d},
in two dimensions we encounter a specific situation that the Hodge star maps 1-forms
on 1-forms again (see (\ref{HodgeE2})).
The latter are, however, in bijection with vector fields (via $\sharp$ and $\flat$).
So, effectively, the Hodge star induces a canonical isomorphism of \emph{vector fields}
on themselves
\begin{equation} \label{vektnaseba}
         \sharp * \flat : \mathfrak X \to \mathfrak X
\end{equation}
A computation gives
\begin{equation} \label{vektnasebakomp}
         \sharp * \flat : (A_x,A_y) \mapsto (-A_y,A_x)
                          \hskip 1cm \text{i.e.} \hskip 1cm
                          \mathbf A \mapsto \mathbf A_{\pi/2}
\end{equation}
We denote the result as $\mathbf A_{\pi/2}$, since it is evident from components
that the new vector is given, at each point, by \emph{rotation} of the vector of the original vector field
by $\pi/2$ in positive sense, i.e. in counter-clock-wise direction.
(This is compatible with the fact that, on 1-forms, the \emph{square} of the star is the \emph{minus} identity,
 see (\ref{HodgeE2nabaze})).

The mapping (\ref{vektnasebakomp}) is equivalent to a handy formula
\begin{equation} \label{vektnasebainac}
         *(\mathbf A \cdot d\mathbf r) =: \mathbf A_{\pi/2} \cdot d\mathbf r
\end{equation}
Then, when combined with (\ref{lavestvorce}), we get the following useful relation
of the gradient field and the Hamiltonian one:
\begin{equation} \label{hamgrad}
         \Ham f = (\grad f)_{\pi/2}
\end{equation}
This is also evident from their expressions in (\ref{grad}) and (\ref{Ham}).
If $f(x,y)$ describes height of a hill somewhere in a landscape,
we know that the vector field $\grad f$ shows, at any point, the direction of the steepest ascent
(and $-\grad f$ the direction of the steepest descent).
Then the vector field $\Ham f$ shows, at any point, the direction of \emph{no} ascent,
i.e. the direction \emph{along contour lines}.

\subsection{Poincar\'e lemma}
\label{poincare}

Recall, first, what this useful lemma states.

The fact $dd=0$ says that if a form $\alpha$ is exact (i.e. if $\alpha = d\beta$),
then it is necessarily closed ($d\alpha =0$).
It turns out that the opposite implication does not, in general, hold.
Poincar\'e lemma, however, guarantees that on a domain which is \emph{contractible to a point}
(for example on a coordinate patch) the opposite implication still does hold.

In 3D vector calculus this leads frequently used and useful statements
(conversions of implications from (\ref{znameidentity}))
\begin{equation} \label{poincarev3D}
          \curl \mathbf a = 0  \ \Rightarrow \ \mathbf a = \grad f
          \hskip 1cm
          \Div \mathbf a = 0  \ \Rightarrow \ \mathbf a = \curl \mathbf b
\end{equation}
When expressed in terms of the diagram (\ref{diagram3D3}), it says that if some arrow
in the bottom line gives zero, the input is necessarily an output of the previous arrow.

And what exactly the deduction gives in 2D vector calculus, i.e. for diagrams (\ref{twoanalogypodrobne})?
This:
\begin{equation} \label{poincarev2D}
          \curl_3 \mathbf a = 0  \ \Rightarrow \ \mathbf a = \grad f
          \hskip 1cm
          \Div \mathbf a = 0  \ \Rightarrow \ \mathbf a = \Ham f
\end{equation}
We can see an application of these facts in Section \ref{2dhydrodynamika}.

\section{2D vector calculus in 2D hydrodynamics}
\label{2dhydrodynamika}

2D vector calculus may be useful, as an example, in description of 2D flows in hydrodynamics.
Such flows are described in terms of \emph{velocity field} $\mathbf v = (v_x,v_y)$.
It results from continuity equation (which encodes conservation of mass under the flow)
that if the mass density $\rho$ is constant, the velocity field is divergence-free:
\begin{equation} \label{divjenula}
          \Div \mathbf v = 0
          \hskip 2cm
          \text{\emph{incompressible} fluid}
\end{equation}

\subsection{Incompressible fluid and stream function}
\label{streamfunction}

Compare the second implications in (\ref{poincarev3D}) and (\ref{poincarev2D}).
We see that if a vector field in 2D happens to be divergence-free,
it is not the curl of an arbitrary vector field, as we are accustomed to in 3D,
but rather it is \emph{Hamiltonian} field generated by an arbitrary function $f(x,y)$.
\begin{equation} \label{hydrodynvhamil2D}
          \Div \mathbf v = 0  \ \Rightarrow \ \mathbf v = \Ham f
          \hskip 1cm \text{i.e.} \hskip 1cm
          (v_x,v_y) = ( - \partial_y f, \partial_x f)
\end{equation}
One can say that the velocity field is \emph{potential} in the sense that components of the field
may be computed in terms of derivatives of some quantity, namely of a function (potential) $f$, here.
The field is, however, not computed as the usual gradient of the potential
(as it is the case, for example, for electric field in electrostatics),
but rather other combinations of partial derivatives are needed.
In particular what we found is that one is to construct Hamiltonian field from $f$.

Now, what is the physical meaning of the function $f$?

First, one easily checks that directional derivative of the function $f$
along \emph{streamlines} of the flow vanishes:
\begin{equation} \label{vsmereprudnic1}
          \dot f \equiv \frac{d}{dt} f(\mathbf r (t))
                 = \dot x (\partial_x f) + \dot y (\partial_y f)
                 = \dot x v_y - \dot y v_x
                 = 0
\end{equation}
since $\mathbf v = (\dot x,\dot y)$ on streamline $\mathbf r (t) \equiv (x(t),y(t))$.

So the function $f$ is \emph{constant} on streamlines. If we were able to find $f$,
we could determine the \emph{shape of streamlines}, i.e. the pattern of our 2D flow
\footnote{The \emph{direction} of the flow as well as the \emph{magnitude of the velocity} of the flow
          cannot be detected in this way.
          The picture shows, at the same time, \emph{contour lines} of the hill with height given by $f(x,y)$,
          mentioned at the end of Section \ref{kolmefield}.
          If the stream function were used in the role of the height of a fictitious hill,
          the fluid would flow, on the corresponding tourist map, along contour lines.}
          , from the equation
\begin{equation} \label{vsmereprudnic3}
                f(x,y) = \ \text{const.}
\end{equation}
That's why $f$ is known as the \emph{stream function}.

Second, we should realize the \emph{value} itself of the function $f$ at a point $A$ cannot carry
any direct physical meaning since it is evident that there is a freedom in additive constant for $f$.
What can carry a real physical meaning, however, is the \emph{difference}
of values at two points $A$ and $B$, i.e. $f(B)-f(A)$.
In electrostatics, where the (electric) field is (minus) the gradient of the potential,
we get \emph{voltage} between the two points in this way.
What we get here, where the velocity field turns out to be Hamiltonian field corresponding to this ``potential''?

It turns out that it gives net \emph{area flux} of the fluid per unit time through a path
connecting the two points.
So if we draw two streamlines which pass through the two points,
it represents the net area flux under a bridge over the ``river'', which runs between the two streamlines.

It may be also seen on a picture (see Fig.\ref{Fig1}) and then computed ``on fingers''.
\begin{figure}[tb]
\begin{center}
\includegraphics[scale=0.40]{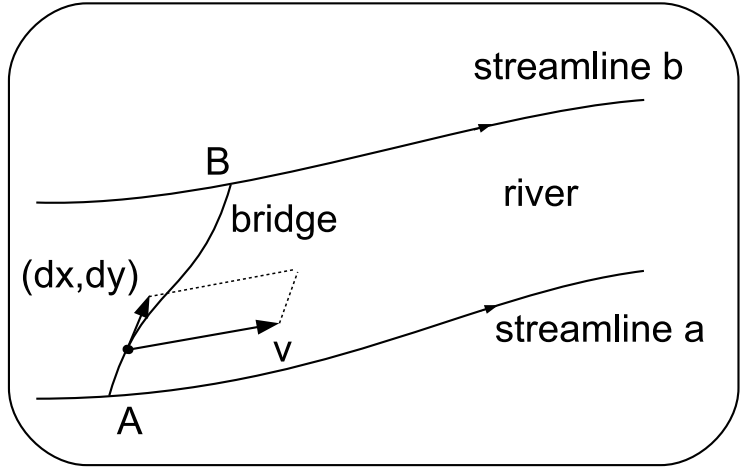}
\caption{Streamline $a$ passes through point $A$, streamline $b$ passes through point $B$.
          Streamlines in between them form a ``river''.
          A bridge over the troubled water connects the point $A$ with the point $B$.
          The area of the rhomboid is the area flux (per unit time) of the water
          underneath the small piece $(dx,dy)$ of the bridge.}
\label{Fig1}
\end{center}
\end{figure}
Indeed, the second identity in (\ref{lavestvorce}) gives
\begin{equation} \label{lavestvorce2}
                 df = - * (\mathbf v \cdot d\mathbf r ) = - * (v_xdx+v_ydy) = v_ydx - v_xdy
\end{equation}
where we used (\ref{HodgeE2nabaze}). So if we step forth (along the path connecting $A$ and $B$)
a small vector $(dx,dy)$, the value of the function $f$ changes by $df = v_ydx - v_xdy$.
Now the expression on the r.h.s. is nothing but the \emph{area of the rhomboid} spanned by edges $(dx,dy)$ and $(v_x,v_y)$.
And this is exactly the area flux of a fluid with velocity $(v_x,v_y)$ through a small piece $(dx,dy)$ of the connecting line,
i.e. the \emph{area}, which flows through the piece per unit time.
(In 3D situation it would be \emph{volume flux}, i.e. the \emph{volume} of the fluid which flows through the small surface $d\mathbf S$
and the corresponding expression were $\mathbf v \cdot d\mathbf S$.)
And the \emph{net change} of the value of $f$ going from $A$ to $B$ is then corresponding \emph{integral}
along the connecting path which gives the \emph{net area flux} for the entire path.

\subsection{Vorticity-free flow and its potential}
\label{potencial}

Now let us turn our attention to the first implication in (\ref{poincarev2D}).
We see that if a vector field $\mathbf a$ in 2D has vanishing $\curl_3 \mathbf a$,
it is necessarily a \emph{gradient} of some function $g(x,y)$. So for the velocity field
$\mathbf v$ it holds
\begin{equation} \label{hydrodynvnevirove2D}
          \curl_3 \mathbf v = 0  \ \Rightarrow \ \mathbf v = \grad g
          \hskip 1cm \text{i.e.} \hskip 1cm
          (v_x,v_y) = ( \partial_x g, \partial_y g)
\end{equation}
But what is the \emph{physical} meaning of the fact $\curl_3 \mathbf v = 0$?

In 3D case, \emph{vector} field $\curl \mathbf v$ is known as the \emph{vorticity} field.
It turns out that, in a given point, it represents twice the vector of the angular velocity,
by which the ``droplet'' centered in the point rotates.

The \emph{scalar} function $\curl_3 \mathbf v$ is the 2D version of the vorticity of the flow.
(It is called vorticity as well.) If we placed a small body made of corc
somewhere on the surface of the 2D flow, it would rotate with angular velocity $2 \ \curl_3 \mathbf v$.

So the condition $\curl_3 \mathbf v = 0$ describes a \emph{vorticity-free} flow.
As we see from (\ref{hydrodynvnevirove2D}), such flow is also ``potential'' one,
moreover in the standard sense, now
(i.e. that the vector field is computed as the gradient of the potential).

\subsection{Where from complex analysis arises}
\label{complexna}

Consider a 2D flow which happens to be ``incompressible''
and \emph{at the same time} vorticity-free.
Then, according to (\ref{hydrodynvhamil2D}) and (\ref{hydrodynvnevirove2D})
it holds
\begin{equation} \label{sucasne}
          \mathbf v = \Ham f = \grad g
          \hskip 1cm \text{i.e.} \hskip 1cm
          (v_x,v_y) = ( - \partial_y f, \partial_x f) = ( \partial_x g, \partial_y g)
\end{equation}
The last equality sign in (\ref{sucasne}) expresses, however, exactly the \emph{Cauchy-Riemann equations};
they say that (complex) function of complex variable
\begin{equation} \label{cauchyriemann}
          h(z)
          \hskip 1cm
          h = f+ig
          \hskip 1cm
          z = x+iy
\end{equation}
is analytic (holomorphic). We can easily check that velocity field is hidden in its derivative
\begin{equation} \label{hprime}
          h'(z) = v_y+iv_x
\end{equation}
Indeed,
\begin{equation} \label{hprimeindeed}
          2\partial_z h(z) = (\partial_x - i\partial_y)(f+ig) =  2(v_y+iv_x)
\end{equation}
One can learn much more on this in standard textbooks on hydrodynamics, e.g. see \cite{LanLif87}.

\section{Conclusion}
\label{zaver}

In this paper the 2D analogue of the standard 3D vector calculus is discussed.
The approach is based on using of differential forms.
The standard logic, which helps so much for understanding 3D vector calculus, is repeated.
The approach is completely ``intrinsic'', it does not use any embedding of the 2D space into
an auxiliary ambient 3D space.
We show what differential operators occur here (analogues of grad, curl and div from 3D),
how they are related and in which integral identities they may be found.
All results are also presented in standard language (i.e. with no forms whatsoever).


\vskip 2cm
The text is the author's translation into English of a paper
\vskip .2cm \noindent
\emph{M.Fecko}:
\newline \noindent
\emph{Vektorov\'a anal\'yza v dvoch rozmeroch},
\newline \noindent
\emph{
Kvaternion, 1-2, 19-30, (2021)}
\vskip .2cm \noindent
originally published in Slovak in journal \emph{Kvaternion}
\newline
(Institute of Mathematics, Faculty of Mechanical Engineering, Brno University of Technology, Czech Republic),
\newline
see
\newline
http://kvaternion.fme.vutbr.cz/index.html
\newline
http://kvaternion.fme.vutbr.cz/2021/kv21\_1-2\_fecko\_web.pdf

\end{document}